\documentclass{article}
\usepackage{PRIMEarxiv}
\usepackage[utf8]{inputenc}
\usepackage[T1]{fontenc}
\usepackage{hyperref}
\usepackage{url}
\usepackage{booktabs}
\usepackage{amsfonts}
\usepackage{amsmath}
\usepackage{amssymb}
\usepackage{nicefrac}
\usepackage[expansion=false]{microtype}
\usepackage{fancyhdr}
\usepackage{graphicx}
\usepackage{enumitem}
\graphicspath{{media/}}

\title{Mathematical Models of Traffic Flow \\ at a Signalized Intersection}

\author{
  Akbota Senkebayeva \\
  Kazakh-British Technical University \\
  Almaty, Kazakhstan \\
  \texttt{akbota.senkebayeva@gmail.com}
}

\begin{document}
\maketitle

\begin{abstract}
This paper presents two one-dimensional mathematical models describing automobile traffic flow on straight road segments at a signalized intersection. When the traffic light is permissive, the flow density and velocity are obtained by solving an initial-boundary value problem for a first-order hyperbolic system. When the signal is prohibitive, the same quantities are governed by a mixed system comprising a second-order parabolic equation for the velocity and a first-order equation for the density.
\end{abstract}

\keywords{traffic flow \and mathematical model \and signalized intersection
\and hyperbolic equations \and parabolic equations \and initial-boundary value
problem}

\section{Introduction}
In this paper two one-dimensional mathematical models
are derived that describe the motion of automobile traffic on straight road
segments under permissive and prohibitive traffic light signals at an
intersection.

Under a permissive traffic light signal, the density and velocity of the
traffic flow are determined as the solution to an initial-boundary value
problem for a system of two first-order hyperbolic equations.

Under a prohibitive traffic light signal, the density and velocity of the
traffic flow are determined as the solution to an initial-boundary value
problem for a system of two equations consisting of a second-order parabolic
equation for the flow velocity and a first-order equation for the flow
density.

In this study, we restrict ourselves to motion only on straight road segments
and only in one direction.

In this formulation, the problem of traffic flow under permissive and
prohibitive traffic light signals at an intersection has not been previously
considered. The concluding part of the paper presents the solution
to one of the posed model problems.

\section{Problem Statement}
The mathematical models proposed in this paper belong
to the group of phenomenological macroscopic models and are based on modified
gas dynamics equations with one spatial variable. The latter have the
form~\cite{ovsiannikov1977}:
\begin{equation}\label{eq:78}
  \rho\!\left(\frac{\partial v}{\partial t}
  + v\frac{\partial v}{\partial x}\right)
  = -\frac{\partial p}{\partial x} + \rho\,F,
\end{equation}
\begin{equation}\label{eq:79}
  \frac{\partial\rho}{\partial t}
  +\frac{\partial}{\partial x}(\rho v)=0, \qquad p=P(\rho).
\end{equation}

Here $\rho(x,t)$ is the density, $v(x,t)$ is the velocity, and $p(x,t)$ is
the pressure of the continuum at a given point $x$ of the continuum at time
$t>0$. The density of external mass forces $F$ and the function $P(s)$ are
assumed to be given. Everywhere below we assume that $F$ is a given positive
function of velocity:
\begin{equation}\label{eq:80}
  F(v)=F_{0}=\mathrm{const}>0 \;\;\text{for}\;\; v<v_{*}-\delta,\qquad
  F(v)=0 \;\;\text{for}\;\; v>v_{*},
\end{equation}
where $v_{*}$ is the maximum permitted speed on the highway outside traffic
control systems.

As a rule, for the system of equations \eqref{eq:78}--\eqref{eq:79}, the
Cauchy problem
\begin{equation}\label{eq:81}
  v(x,0)=v_{0}(x), \qquad \rho(x,0)=\rho_{0}(x)
\end{equation}
is considered on the entire interval $-\infty<x<\infty$ for $t>0$, or the
initial-boundary value problem
\begin{equation}\label{eq:82}
  v(x,0)=v_{0}(x), \qquad \rho(x,0)=\rho_{0}(x), \qquad x>0,
\end{equation}
\begin{equation}\label{eq:83}
  v(0,t)=v^{0}(t), \qquad \rho(0,t)=\rho^{0}(t), \qquad t>0
\end{equation}
on the half-line $x>0$ for $t>0$ under the condition $v_{0}(x) \geq 0$ and
$v^{0}(t) \geq 0$.

When the prohibitive (red) traffic light signal is activated at the
intersection $x=x_{0}>0$, the traffic flow must stop at that point at a
certain time $t=t_{0}$. Moreover, it must begin braking at the flashing green
light at time $t=t_{0}-\tau_{0}$, slightly before the established stopping
time $t=t_{0}$.

Controlling this process in the mathematical model is only possible with an
additional boundary condition. Namely, a value $h$ and a boundary
\[
  \Gamma = \{x=\gamma(t),\; t_{0}-\tau_{0} < t < t_{0}\}, \qquad
  \gamma(t_{0}-\tau_{0}) = x_{0}-h, \qquad
  \gamma(t) = x_{0} \;\;\text{for}\;\; t \geq t_{0},
\]
are chosen, on which
\begin{equation}\label{eq:84}
  v\!\left(\gamma(t),\,t\right) = V(t) \qquad \text{for}\;\; t \geq t_{0},
\end{equation}
where $V(t)$ is a continuous function such that
$V(t_{0}-\tau_{0}) = v(x_{0}-h,\,t_{0}-\tau_{0}-0)$ and $V(t)=0$ for
$t \geq t_{0}$ (at time $t=t_{0}$ the traffic flow has stopped before the traffic
light and remains at rest until the permissive signal is activated).

Unfortunately, the system of equations \eqref{eq:78}, \eqref{eq:79} is
hyperbolic and does not admit any additional boundary conditions. Therefore,
for $t>t_{0}-\tau_{0}$, the motion of the traffic flow must be described by
different differential equations.

For example, for the initial-boundary value problem
\eqref{eq:78}, \eqref{eq:79}, \eqref{eq:82}, \eqref{eq:83}, we consider the
viscous gas equations
\begin{equation}\label{eq:85}
  \rho\!\left(\frac{\partial v}{\partial t}
  + v\frac{\partial v}{\partial x}\right)
  = \mu\,\frac{\partial^{2} v}{\partial x^{2}}
  - \frac{\partial p}{\partial x},
\end{equation}
\begin{equation}\label{eq:86}
  \frac{\partial\rho}{\partial t}
  + \frac{\partial}{\partial x}(\rho v) = 0, \qquad p = P(\rho),
\end{equation}
supplemented with the following boundary and initial conditions:
\begin{equation}\label{eq:87}
  v(0,t) = v^{0}(t), \qquad \rho(0,t) = \rho^{0}(t), \qquad
  t > t_{0}-\tau_{0},
\end{equation}
\begin{equation}\label{eq:88}
  v\!\left(\gamma(t),\,t\right) = V(t), \qquad t > t_{0}-\tau_{0},
\end{equation}
\begin{equation}\label{eq:89}
  \begin{aligned}
    &v(x,\,t_{0}-\tau_{0}+0) = v(x,\,t_{0}-\tau_{0}-0),\\
    &\rho(x,\,t_{0}-\tau_{0}+0) = \rho(x,\,t_{0}-\tau_{0}-0),
    \qquad 0 < x < x_{0}-h,
  \end{aligned}
\end{equation}
under the condition that the given function $V$ satisfies
\begin{equation}\label{eq:90}
  V(t_{0}-\tau_{0}) = v(x_{0},\,t_{0}-\tau_{0}+0)
  \qquad\text{and}\qquad
  V(t_{0}) = 0.
\end{equation}

In \eqref{eq:85}--\eqref{eq:89}, $\mu=\mathrm{const}>0$ is the viscosity of
the traffic flow, and $\{v(x,t_{0}-\tau_{0}-0),\,\rho(x,t_{0}-\tau_{0}-0)\}$
is the limit of the solutions $\{v(x,t),\,\rho(x,t)\}$ of the
initial-boundary value problem \eqref{eq:78}, \eqref{eq:79}, \eqref{eq:82},
\eqref{eq:83} on the interval $(0,\,t_{0}-\tau_{0})$ as
$t\nearrow t_{0}-\tau_{0}$. The question of whether to disable the vehicle
acceleration ($F=0$) on the entire interval or only near the traffic light
remains open. In the proposed equation \eqref{eq:85} we assume $F=0$ on the
entire interval of motion.

Problem \eqref{eq:85}--\eqref{eq:89} is solved up to the moment
$t=t_{0}+\tau_{1}$ when the permissive green traffic light signal is
activated.

The second traffic flow on the interval $x_{0}-h < x < \infty$ for
$t_{0}-\tau_{0} < t < t_{0}+\tau_{1}$ continues its motion, since all
vehicles on this interval manage to pass the point $x=x_{0}$ before the
prohibitive red signal is activated. In this case, the flow motion is
described by the original system \eqref{eq:78}, \eqref{eq:79} and initial
condition \eqref{eq:82} for $x_{0}-h < x < \infty$, and the boundary
condition \eqref{eq:83} must be replaced by the following condition:
\begin{equation}\label{eq:91}
  v(x_{0}-h,\,t) = 0, \qquad \rho(x_{0}-h,\,t) = 0, \qquad
  t_{0}-\tau_{0} < t < t_{0}+\tau_{1},
\end{equation}
expressing the absence of vehicle flow through the boundary $x=x_{0}-h$ for
the second flow.

When the permissive green traffic light signal is activated at
$t=t_{0}+\tau_{1}$, the procedure repeats, but with new initial conditions
\begin{equation}\label{eq:92}
  \begin{aligned}
    &v(x,\,t_{0}+\tau_{1}+0) = v(x,\,t_{0}+\tau_{1}-0),\\
    &\rho(x,\,t_{0}+\tau_{1}+0) = \rho(x,\,t_{0}+\tau_{1}-0),
    \qquad 0 < x < \infty,
  \end{aligned}
\end{equation}
where $\{v(x,t_{0}+\tau_{1}-0),\,\rho(x,t_{0}+\tau_{1}-0)\}$ for
$0 < x < x_{0}$ is the limit of the solutions $\{v(x,t),\,\rho(x,t)\}$ of
the initial-boundary value problem \eqref{eq:85}--\eqref{eq:89}, and for
$x_{0} < x < \infty$ is the limit of the solutions $\{v(x,t),\,\rho(x,t)\}$
of the initial-boundary value problem
\eqref{eq:78}--\eqref{eq:80}, \eqref{eq:82}, \eqref{eq:91} on the interval
$t_{0}-\tau_{0} < t < t_{0}+\tau_{1}$ as $t\nearrow t_{0}+\tau_{1}$.

Note that compared to a gas or liquid, the traffic flow has its own specific
features. If a liquid or gas is at rest at the initial moment $t=0$
($v(x,0)=0$), but the initial density is distributed non-uniformly
($\rho(x,0)\neq\mathrm{const}$), then even in the absence of external mass
forces ($F=0$), the liquid or gas will begin to spread for $t>0$. Cars,
however, regardless of how unevenly they are arranged, will remain stationary
without an external accelerating force $F$. The factor that causes a liquid or
gas to spread (equalize density) is the pressure gradient in equations
\eqref{eq:78} or \eqref{eq:85}. The only way to avoid such behavior of a
continuum is the postulate
\begin{equation}\label{eq:93}
  P(s) \equiv \mathrm{const}.
\end{equation}

The initial-boundary value problem
\eqref{eq:78}--\eqref{eq:80}, \eqref{eq:82}, \eqref{eq:83}, \eqref{eq:93},
describing the motion of the traffic flow on the interval $(0,\,t_{0}-\tau_{0})$,
will be called \textbf{Problem~$A_{1}$}.

The initial-boundary value problem \eqref{eq:85}--\eqref{eq:89},
\eqref{eq:93}, describing the motion of the traffic flow on the interval
$(t_{0}-\tau_{0},\,t_{0}+\tau_{1})$, will be called \textbf{Problem~$B_{1}$}.

Finally, the initial-boundary value problem
\eqref{eq:78}--\eqref{eq:80}, \eqref{eq:82}, \eqref{eq:91}, \eqref{eq:93},
describing the motion of the traffic flow on the interval
$(t_{0}-\tau_{0},\,t_{0}+\tau_{1})$, will be called \textbf{Problem~$C_{1}$}.

Problems $A_{1}$, $B_{1}$, and $C_{1}$ form the \textbf{first mathematical
model} of traffic flow on a straight road segment.

In the first model, the differential equations of continuum motion depend on
the traffic light signal at the point $x=x_{0}$. Before this signal, the
traffic flow is described by a hyperbolic system, and after the signal, part
of the flow is described by a mixed-type system. Another approach is also
possible, where the motion of the traffic flow is described by the same system
of differential equations throughout. Namely, instead of equation \eqref{eq:78}
for $0<t<t_{0}-\tau_{0}$ on the segment $0<x<x_{0}-h$ and for
$t_{0}-\tau_{0}<t<t_{0}+\tau_{1}$ on the segment $x_{0}-h<x<\infty$, we
consider the equations
\begin{equation}\label{eq:94}
  \rho\!\left(\frac{\partial v}{\partial t}
  + v\frac{\partial v}{\partial x}\right)
  = \mu\,\frac{\partial^{2} v}{\partial x^{2}} + \rho\,F,
\end{equation}
\begin{equation}\label{eq:95}
  \frac{\partial\rho}{\partial t}
  + \frac{\partial}{\partial x}(\rho v) = 0,
\end{equation}
keeping the same boundary and initial conditions \eqref{eq:82}, \eqref{eq:83}
and \eqref{eq:91} (\textbf{Problems~$A_{2}$ and $C_{2}$} respectively).

For $t_{0}-\tau_{0}<t<t_{0}+\tau_{1}$ in the domain
$\{t_{0}-\tau_{0}<t<t_{0}+\tau_{1},\;0<x<\gamma(t)\}$, the
initial-boundary value problem \eqref{eq:87}--\eqref{eq:89} is solved for
the equations
\begin{equation}\label{eq:96}
  \rho\!\left(\frac{\partial v}{\partial t}
  + v\frac{\partial v}{\partial x}\right)
  = \mu\,\frac{\partial^{2} v}{\partial x^{2}} + \rho\,F,
\end{equation}
\begin{equation}\label{eq:97}
  \frac{\partial\rho}{\partial t}
  + \frac{\partial}{\partial x}(\rho v) = 0,
\end{equation}
which coincide with equations \eqref{eq:94} and \eqref{eq:95} if we set
$F=0$ in them (\textbf{Problem~$B_{2}$}).

Problems $A_{2}$, $B_{2}$, and $C_{2}$ form the \textbf{second mathematical
model} of traffic flow on a straight road segment.

\section{Generalized Solutions}
The problem under study is extremely complex from a mathematical point of
view (existence and uniqueness of solution), and even in the simplest
situations one should not expect the solution to have all the derivatives
that form the differential equations. Therefore, it is quite reasonable to
study weak (generalized) solutions, where the differential equations take the
form of conservation laws.

\section{Solution of the Initial-Boundary Value Problem Describing
Traffic Flow on a Straight Road Segment}
Here, one solution of the initial-boundary value problem
\eqref{eq:78}--\eqref{eq:83} is proposed. In this problem, we assume that
\begin{equation}\label{eq:98}
  \rho(x,0) = \rho_{0}(x) \geq 0, \qquad
  v(x,0) = v_{0}(x) \geq 0 \qquad (x>0),
\end{equation}
\begin{equation}\label{eq:99}
  \rho(0,t) = \rho^{0}(t) \geq 0, \qquad
  v(0,t) = v^{0}(t) \geq 0 \qquad (t>0).
\end{equation}

Consider new independent variables $(\xi,\,t)$ introduced by the formula
\begin{equation}\label{eq:100}
  \xi = \int_{0}^{x}\rho(s,t)\,ds, \qquad t = t.
\end{equation}

Let the new unknown functions
\begin{equation}\label{eq:101}
  \hat{\rho}(\xi,t) = \rho(x,t), \qquad \hat{v}(\xi,t) = v(x,t)
\end{equation}
satisfy the following system of equations:
\begin{equation}\label{eq:102}
  \frac{\partial\hat{v}}{\partial t}
  + a(t)\frac{\partial\hat{v}}{\partial\xi} = F,
\end{equation}
\begin{equation}\label{eq:103}
  \frac{\partial\hat{\rho}}{\partial t}
  + a(t)\frac{\partial\hat{\rho}}{\partial\xi}
  + \hat{\rho}^{2}\frac{\partial\hat{v}}{\partial\xi} = 0,
\end{equation}
where $a(t) = \rho^{0}(t)\cdot v^{0}(t)$.

In deriving \eqref{eq:102}, \eqref{eq:103}, we used the formulas:
\begin{align*}
  \frac{\partial\rho}{\partial t}
  &= \frac{\partial\hat{\rho}}{\partial t}
     + \frac{\partial\hat{\rho}}{\partial\xi}\cdot\frac{\partial\xi}{\partial t}
   = \frac{\partial\hat{\rho}}{\partial t}
     + \frac{\partial\hat{\rho}}{\partial\xi}
       \left(\int_{0}^{x}\frac{\partial\rho}{\partial t}(s,t)\,ds\right) \\
  &= \frac{\partial\hat{\rho}}{\partial t}
     - \frac{\partial\hat{\rho}}{\partial\xi}
       \left(\int_{0}^{x}\frac{\partial}{\partial s}(\rho v)(s,t)\,ds\right)
   = \frac{\partial\hat{\rho}}{\partial t}
     - \frac{\partial\hat{\rho}}{\partial\xi}\hat{\rho}\hat{v}
     + \frac{\partial\hat{\rho}}{\partial\xi}\cdot a(t)
\end{align*}
and
\[
  \frac{\partial\rho}{\partial x}
    = \frac{\partial\hat{\rho}}{\partial\xi}\cdot\frac{\partial\xi}{\partial x}
    = \hat{\rho}\,\frac{\partial\hat{\rho}}{\partial\xi}; \qquad
  \frac{\partial v}{\partial x}
    = \hat{\rho}\,\frac{\partial\hat{v}}{\partial\xi}; \qquad
  \frac{\partial v}{\partial t}
    = \frac{\partial\hat{v}}{\partial t}
    - \frac{\partial\hat{v}}{\partial\xi}\hat{\rho}\hat{v}
    + \frac{\partial\hat{v}}{\partial\xi}\,a(t).
\]

We will seek solutions in which $\rho(x,t)\geq 0$. Therefore, the
transformation \eqref{eq:100} maps the half-line $\{x>0\}$ to the half-line
$\{\xi>0\}$. On the boundary $\{\xi=0\}$, the boundary conditions
\begin{equation}\label{eq:104}
  \hat{\rho}(0,t) = \rho^{0}(t), \qquad \hat{v}(0,t) = v^{0}(t)
\end{equation}
hold, which follow from conditions \eqref{eq:99}. Since
$\hat{\rho}(\xi,0)=\hat{\rho}_{0}(\xi)=\rho(x,0)=\rho_{0}(x)$ and
$\hat{v}(\xi,0)=\hat{v}_{0}(\xi)=v(x,0)=v_{0}(x)$, to determine the
initial values $\hat{\rho}_{0}(\xi)$ and $\hat{v}_{0}(\xi)$ we must find
the transformation \eqref{eq:100} at $t=0$:
$\xi=\int_{0}^{x}\rho_{0}(s)\,ds$, and compute its inverse
$x=\gamma_{0}(\xi)$. The condition for the existence of the inverse
transformation is
\begin{equation}\label{eq:105}
  \rho_{0}(x) > 0 \qquad (x>0).
\end{equation}

Thus, we can define $\hat{\rho}_{0}(\xi)=\rho_{0}(\gamma_{0}(\xi))$ and
$\hat{v}_{0}(\xi)=v_{0}(\gamma_{0}(\xi))$, and set the initial conditions
$\hat{\rho}(\xi,0)=\hat{\rho}_{0}(\xi)$,
$\hat{v}(\xi,0)=\hat{v}_{0}(\xi)$ with known functions $\hat{\rho}_{0}$
and $\hat{v}_{0}$.

\bibliographystyle{unsrt}
\bibliography{references}

\end{document}